\documentclass[11pt,reqno]{amsart}
\usepackage{amsmath,amssymb,amsthm,upref,graphicx,mathrsfs}

\usepackage{color,xcolor}
\usepackage[
  colorlinks=true,
  linkcolor=blue,
  citecolor=blue,
  urlcolor=blue]{hyperref}

\newtheorem{thm}{Theorem}[section]
\newtheorem{cor}[thm]{Corollary}
\newtheorem{lem}[thm]{Lemma}

\newtheorem{exaple}[thm]{Example}
\numberwithin{equation}{section}

\begin{document}

\leftline{ \scriptsize}

\vspace{1.3 cm}
\title
{A remark on an equivalence of two versions of Polya's permanent problem }
\author{Ratsiri Sanguanwong $^{a}$ and Kijti Rodtes$^{b,\ast}$ }
\thanks{{\scriptsize
		$^{\ast}$ Corresponding author.\\
		E-mail addresses: ratsiris60@email.nu.ac.th (Ratsiri Sanguanwong), kijtir@nu.ac.th (Kijti Rodtes).\\
		$^{A}$ Department of Mathematics, Faculty of Science, Naresuan University, Phitsanulok 65000, Thailand.\\}}
\hskip -0.4 true cm

\maketitle


\begin{abstract}In 2004, some equivalent versions of Polya's permanent problem were listed in 24 versions. However, there is a flaw on the theorem that affirms an equivalence of version 11 and 12.  In order to correct the slip, we provide a characterization when the determinant and permanent of a given nonnegative matrix are equal. Moreover, the new characterization yields necessary and sufficient conditions for the equality of the determinant and permanent of powers of a given nonnegative matrix.
\end{abstract}

\vskip 0.2 true cm


\pagestyle{myheadings}
\markboth{\rightline {\scriptsize Ratsiri Sanguanwong and  Kijti Rodtes}}
{\leftline{\scriptsize }}
\bigskip
\bigskip


\vskip 0.4 true cm

\section{Introduction}

The Polya's  permanent problem has been studied intensively since 1913. The main idea of the problem is to calculate the permanent of a given matrix via the determinant of the transformed matrix. The following problems are related to Polya's work \cite{GP1913}.

\begin{quote}
	\textbf{Problem 1} For a given $n \times n$ matrix $A=(a_{ij})$, does there exist a uniform way of affixing $\pm$ sign such that $\operatorname{per}(a_{ij})=\det(\pm a_{ij})$?\\
	\textbf{Problem 2} Does there exist a transformation $\phi:M_n(\mathbb{F})\rightarrow M_m(\mathbb{F})$ satisfying $\operatorname{per}(A)=\det(\phi(A))$ for all $A \in M_n(\mathbb{F})$?\\
	\textbf{Problem 3} For a given $n \times n$ $\{0,1\}-$matrix $A$, does there exist a matrix $B$ obtained by changing some of the $1$ entries of $A$ to $-1$ such that $\operatorname{per}(A)=\det(B)$?
\end{quote}

In \cite{GZ1913}, Szego observed that, if $n \neq 2$, the uniform way as in Problem 1 does not exist. In 1961, Problem 2 was also proved negative for $m=n>2$ by Marcus and Minc (cf. \cite{MM1961}).

The motivation of this paper is Problem 3. In \cite{WM2004}, McCuaig showed that Problem 3 has other 24 versions while itself is referred as Version 10. Most of these versions are in the form of graph theoretical questions. Our interest is in the following two versions.

\begin{quote}
	\textbf{Version 11} Does a given digraph have a dicycle of even length?\\
	\textbf{Version 12} For a given square nonnegative matrix $A$, does $\operatorname{per}(A)=\det(A)$?
\end{quote}

These two versions are claimed to be equivalent to each other by the statements a) and b) in Theorem 10 \cite{WM2004} which stated as follow.

\begin{thm}
	Let $A$ be a square nonnegative matrix such that $\det(A)$ has a positive term. Then the bipartite graph $G$ of $A$ has a perfect matching. Let $D$ be a bipartite graph corresponding to $G$. The following statements are equivalent.
	\begin{itemize}
		\item[a)] All dicycles of $D$ have odd length.
		\item[b)] $\operatorname{per}(A)=\det(A)$.
	\end{itemize}
\end{thm}

The digraph $D$ in this theorem is the same graph as $G(A)$ defined in this paper (see Section 3). However, if we choose $\sigma = (1234)(5678)$, we can see that $\operatorname{per}(P_\sigma)=1=\det(P_\sigma)$ while $G(P_\sigma)$ contains 2 cycles of length 4.  Since Theorem 10 \cite{WM2004} contains neither a proof nor a citation, we have no clue about a cause of the slip. Thus, in this paper, we assume one more assumption to make this theorem holds.


In this paper, for every nonnegative matrix $A$ with $\operatorname{per}(A)\neq 0$, a necessary and sufficient condition for $\operatorname{per}(A)=\lvert \det(A) \rvert$ is investigated. We also give a sufficient condition for $\operatorname{per}(A)=\det(A)$ so that the condition does not deal with a graph. Finally, by using those results, a characterization for $\operatorname{per}(A^k) = \det(A^k)$ for all $k \in \mathbb{N}$ is also established. 

\section{Preliminaries}

In this section, basic definitions, notations, and properties are given. Firstly, we begin with those of matrices.

For any matrix $A$, we denote $A_{ij}$ the $(i,j)$-entry of $A$. For a given set $S$, a matrix $A$ is called an \textit{$S$-matrix} if each entry of $A$ is an element of $S$. An $S$-matrix $A$ is said to be a \textit{sign pattern matrix} if $S=\{-,0,+\}$. A matrix $A$ is said to be \textit{complex} (resp. \textit{real},\textit{ nonnegative}) if $S$ is the set of all complex (resp. real, nonnegative) numbers. For a sign pattern matrix $\Gamma$, the set $Q(\Gamma)$ denotes the set of real matrices whose each entry has the same sign as $\Gamma$. 

The \textit{permanent} and the \textit{determinant} of an $n \times n$ complex matrix $A$, denoted by $\operatorname{per}(A)$ and $\det(A)$, respectively, are defined as
\begin{equation*}
\hbox{$\operatorname{per}(A)=\sum\limits_{\sigma\in S_n}\prod\limits_{i=1}^n A_{i\sigma(i)}$ and $\det(A)=\sum\limits_{\sigma\in S_n}\operatorname{sgn}(\sigma)\prod\limits_{i=1}^n A_{i\sigma(i)}$.}
\end{equation*}
For any $n\times n$ matrices $A,B$, we have $\det(BA)=\det(AB)=\det(A)\det(B)$. However, this property does not hold for the permanent. For convenience, for each positive integer $n$ and each $\sigma\in S_n$, we denote $[n]:=\{1,\dots,n\}$ and $L_\sigma(A):=\prod\limits_{i=1}^nA_{i \sigma(i)}$. Recall that for an $n \times n$ matrix $A$,
$$\det(A)=\sum\limits_{\sigma\in A_n}L_\sigma(A) - \sum\limits_{\sigma\in S_n\setminus A_n}L_\sigma(A)$$
and
$$\operatorname{per}(A)=\sum\limits_{\sigma\in A_n}L_\sigma(A) + \sum\limits_{\sigma\in S_n\setminus A_n}L_\sigma(A).$$

So, we can see that $\operatorname{per}(A) = \det(A)$ if and only if $\sum\limits_{\sigma\in S_n\setminus A_n}L_\sigma(A)=0$. Similarly, $\operatorname{per}(A) = -\det(A)$ if and only if $\sum\limits_{\sigma\in A_n}L_\sigma(A)=0$.

Let $A$ be a square nonnegative matrix. We have that $\operatorname{per}(A)\geq \lvert \det(A) \rvert$. Thus, if $\operatorname{per}(A)=0$, then $\det(A)=0$. The characterization for $\operatorname{per}(A)=0$ can be obtained by Frobenius-K\"{o}nig's Theorem; that is, $\operatorname{per}(A)=0$ if and only if $A$ contains $r \times s$ zero submatrix with $r+s=n+1$. In this paper, we only work with nonnegative matrices whose permanent is not zero.

A \textit{permutation matrix} is a square $\{0,1\}$-matrix whose each row and each column contains exactly one $1$. Let $\mathbb{P}_n$ be the set of all $n \times n$ permutation matrices. We have that $\mathbb{P}_n$ is a group under multiplication and isomorphic to the symmetric group $S_n$. For $\sigma \in S_n$, we denote $P_\sigma$ the permutation matrix obtained by permuting each row $i$ to $\sigma(i)$. This means that
\begin{equation*}\label{pm}
(P_\sigma)_{ij} = \left\{
\begin{array}{ll}
1, & \hbox{if $i=\sigma(j)$,} \\
0, & \hbox{otherwise.}
\end{array}
\right.
\end{equation*}
Thus $\operatorname{per}(P_\sigma) = \prod\limits_{i=1}^n(P_\sigma)_{i\,\sigma^{-1}(i)}$. Because of this fact, we should be careful of a confusion of subscripts in each definition. For any two permutation $\sigma,\tau \in S_n$ and any $n \times n$ matrix $A$, $(P_\sigma A P_\tau)_{ij} = A_{\sigma^{-1}(i)\tau(j)}$ for all $i,j \in [n]$.

Next, some definitions and properties in graph theory are given. In this paper, the results is regarded to a special type of graphs, which are called digraphs. Let $G$ be a graph with a vertex set $V$ and an edge set $E$. A \textit{walk} of length $k$ in $G$ is a finite sequence of vertices and edges $v_1e_1v_2e_2v_3\dots v_{k}e_kv_{k+1}$ such that $v_i$ and $v_{i+1}$ incident with $e_i$ for each $i\in[k]$. If a graph $G$ is simple, we can omit $e_1,\dots,e_k$ from the sequence. Here, $v_1$ is called the \textit{origin vertex} while $v_{k+1}$ is the \textit{terminus vertex}. If $v_1,\dots,v_{k+1}$ are all distinct, a walk $v_1\dots v_{k+1}$ is said to be a \textit{path}. A path $v_1\dots v_{k+1}$, when $k\geq3$, is said to be a \textit{cycle} if $v_{k+1}=v_1$; i.e., it is a path whose the origin vertex and the terminus vertex are the same vertex. 

A \textit{digraph} or \textit{directed graph} is a graph in which each edge is assigned a direction. An edge with a direction is called an \textit{arc}. A cycle in a digraph is called a \textit{dicycle}. Let $G$ be a digraph with a vetex set $\{v_1,\dots,v_n\}$. The \textit{adjacency matrix} of $G$ is an $n \times n$ $\{0,1\}$-matrix $A$ such that $A_{ij}=1$ if and only if $v_iv_j$ is an arc of $G$. So, in case $G$ is simple, all diagonal entries of $A$ are 0. The following theorem discusses a relation between the adjacency matrix of a graph and the number of walks.
\begin{thm}[Theorem 2.2 \cite{CLZB}]
	Let $G$ be a graph with a vertex set $\{v_1, \dots, v_n\}$ and $A$ is the adjacency matrix of $G$. For each positive integer $k$, the number of walks of length $k$ from $v_i$ to $v_j$ in $G$ is the $(i,j)$-entry of $A$.
\end{thm}

For a matrix $A$, the support of $A$ is the set $\{(i,j) \mid A_{ij}\neq 0\}$. Let $A,B$ be $n \times n$ nonnegative matrices. If $\operatorname{per}(A) = \lvert \det(A) \rvert$ and the support of $B$ is a subset of the support of $A$, then $\operatorname{per}(B) = \lvert \det(B) \rvert$. Define a function $\phi$ from $M_n(\mathbb{C})$ onto the set of all $n\times n$ $\{0,1\}$-matrices by $(\phi(A))_{ij}=0$ if and only if $A_{ij}=0$. Let $\Gamma$ be a sign pattern matrix. We can see that $\phi(A) = \phi(B)$ for any matrices $A,B\in Q(\Gamma)$. We denote the $\{0,1\}$-matrix $\phi(A)$ by $\phi(\Gamma)$. Then the following Lemma holds.

\begin{lem}\label{cor1}
	Let $\Gamma$ be a sign pattern matrix. Then the following statements are equivalent:
	\begin{enumerate}
		\item $\operatorname{per}(A) = \lvert \det(A) \rvert$ for all $A \in Q(\Gamma)$,
		\item $\operatorname{per}(A) = \lvert \det(A) \rvert$ for some $A \in Q(\Gamma)$,
		\item $\operatorname{per}(\phi(\Gamma)) = \lvert \det(\phi(\Gamma)) \rvert$.
	\end{enumerate}
\end{lem}

Lemma \ref{cor1} can be extended to the version of complex matrices. Let $\Gamma$ be a $\{0,\ast\}$-matrix. We define the set $Q(\Gamma)$ analogue to the version of sign pattern matrix by $Q(\Gamma)=\{A\in M_n(\mathbb{C})\mid A_{ij}=0\text{ if and only if }\Gamma_{ij}=0\}$. The following lemma is also obtained immediately.

\begin{lem}\label{cor22}
	Let $\Gamma$ be a $\{0,\ast\}$-matrix. Then the following statements are equivalent:
	\begin{enumerate}
		\item $\operatorname{per}(A) = \pm\det(A)$ for all $A \in Q(\Gamma)$,
		\item $\operatorname{per}(A) = \pm\det(A)$ for some $A \in Q(\Gamma)$,
		\item $\operatorname{per}(\phi(\Gamma)) = \pm\det(\phi(\Gamma))$.
	\end{enumerate}
\end{lem}

\section{Main Results}

If we have an $n \times n$ matrix $A$, we can construct the digraph $G(A)$ of $A$ with a vertex set $[n]$ and an edge set $\{ij \mid A_{ij}\neq 0 \text{ and } i\neq j\}$.

\begin{thm}\label{lem1}
	Let $A$ be a nonnegative matrix such that $A_{ii}>0$ for all $i\in[n]$. Then $\operatorname{per}(A)=\det(A)$ if and only if the following statements hold:
	\begin{enumerate}
		\item For all distinct $i,j\in[n]$, $A_{ij}=0$ or $A_{ji}=0$.
		\item The digraph $G(A)$ contains no cycle of even length.
	\end{enumerate}
\end{thm}
\begin{proof}
	Assume that $\operatorname{per}(A)=\det(A)$. Then, for each $\sigma \in S_n,$ we have $L_\sigma(A)\neq 0$ implies $\sigma\in A_n$. Suppose that there exists distinct $i,j\in[n]$ such that neither $A_{ij}$ nor $A_{ji}$ is $0$. Because all diagonal entries are not zero, we obtain that $L_{(ij)}(A)\neq 0$. This contradicts with $\operatorname{per}(A)=\det(A)$ since $(ij)$ is not in $A_n$. Thus (1) is true. To prove (2), suppose that $v_1 v_2 \dots v_k v_1$ is a cycle of length $k$ in $G(A)$. Let $c=(v_1 v_2 \dots v_k) \in S_n$. By the definition of $G(A)$, we have $A_{v_k v_1}\neq 0$ and $A_{v_i v_{i+1}}\neq 0$ for all $i\in[k-1]$. Because all diagonal entries of $A$ are not zero, $L_c(A) \neq0$. So, $c \in A_n$; i.e., $k$ is odd. Therefore (2) holds.
	
	Conversely, suppose that (1) and (2) occur. Let $c$ be a cycle of length $k$ in $S_n$ such that $L_c(A)\neq 0$. Then $c = (i\, c(i)\, \dots \,c^{k-1}(i))$ for some $i \in [n]$. Claim that $k$ is odd. Because $L_c(A)\neq 0$, we can conclude that $$A_{i\, c(i)}, A_{c(i)\,c^2(i)}, \dots, A_{c^{k-2}(i)\, c^{k-1}(i)}, A_{c^{k-1}(i)\, i}$$ are all positive. If $k=2$, then $A_{i \, c(i)}$ and $A_{c(i) \, i}$ are both positive. This contradicts to (1). Thus $k\geq 3$. This means that $i \, c(i) \, \dots \, c^{k-1}(i) \, i$ is a cycle in $G(A)$. Because of (2), $k$ is odd. Let $\tau\in S_n$ such that $L_\tau(A) \neq 0$. If $\tau = c_1 \cdots c_m$ is the disjoint cycles decomposition, then $L_{c_i}(A) \neq 0$ for all $i \in [m]$. Thus each $c_i$ is a cycle of odd length. This implies that $\tau\in A_n$. Since $\tau$ is arbitrary, $\operatorname{per}(A)=\det(A)$.
\end{proof}

By Theorem \ref{lem1}, the following corollary is obtained immediately.

\begin{cor}\label{corthm1}
	Let $A$ be a nonnegative matrix. Then $0\neq \operatorname{per}(A)=\lvert \det(A) \rvert$ if and only if the following statements hold:
	\begin{enumerate}
		\item There exists $\tau \in S_n$ such that $L_\tau(A)\neq 0$. 
		\item For all distinct $i,j\in[n]$, $(AP_\tau)_{ij}=0$ or $(AP_\tau)_{ji}=0$.
		\item The digraph $G(AP_\tau)$ contains no cycle of even length.
	\end{enumerate}
	In particular, if $A$ satisfies $(1)$-$(3)$, then $\operatorname{per}(A)=\operatorname{sgn}(\tau)\det(A)$.
\end{cor}
\begin{proof}
	Suppose that $0\neq \operatorname{per}(A) = \lvert\det(A)\rvert$. The statement (1) is obvious because $\operatorname{per}(A)\neq 0$. Note that $A_{i \tau(i)} = (AP_\tau)_{ii}$ for all $i\in[n]$. This means that all diagonal entries of $AP_\tau$ are positive. Note that $\det(AP_\tau) = \operatorname{sgn}(\tau)\det(A)$. Since $\operatorname{per}(A)=\lvert \det(A) \rvert$, we can conclude that $\det(A)=\operatorname{per}(A)$ if and only if $\tau \in A_n$. This implies that $$\operatorname{per}(AP_\tau)=\operatorname{per}(A)=\lvert \det(A) \rvert = \operatorname{sgn}(\tau)\det(A) = \det(AP_\tau).$$
	By Theorem \ref{lem1}, (2) and (3) occur.
	
	Conversely, suppose that (1)-(3) happen. By Theorem \ref{lem1}, we have $$\operatorname{per}(A)=\operatorname{per}(AP_\tau)=\det(AP_\tau)=\operatorname{sgn}(\tau)\det(A).$$
	Since $\operatorname{per}(A)$ is positive, $\operatorname{per}(A) = \lvert \det(A) \rvert$.
\end{proof}

Note that the number of zero entries in $A$ is the same number as in $AP_\tau$. Therefore, the following corollary is also true.

\begin{cor}\label{corthm2}
	Let $A$ be nonnegative matrix such that $0 \neq \operatorname{per}(A)=\lvert \det(A)\rvert$. Then $A$ contains at least $\frac{n^2-n}{2}$ zero entries.
\end{cor}
The condition 1 of Theorem \ref{lem1} also allows us to conclude that:
\begin{cor}
	Let $A$ be an $n \times n$ positive semidefinite matrix. Then $\operatorname{per}(A)$=$\det(A)$ if and only if $A$ contains a zero rows or $A$ is diagonal.
\end{cor}
\begin{proof}
	Note that every principal submatrix of a positive semidefinite matrix is also positive semidefinite. Suppose that there exists $i\in[n]$ such that $A_{ii}=0$. For each $j\in[n]$ such that $j\neq i$, we have $$\begin{pmatrix}
	A_{ii} & A_{ij}\\
	A_{ji} & A_{jj}
	\end{pmatrix}$$ is positive semidefinite. This implies that the determinant of this matrix is nonnegative. Since $A_{ii}=0$ and $A$ is Hermitian, we have $A_{ij}=0=A_{ji}$. Then all entries in the $i$-th row of $A$ are zero (cf.\cite{FZhangB}). Suppose that $A_{ii}\neq 0$ for all $i\in[n]$. By Theorem \ref{lem1}, for any distinct $i,j\in[n]$, $A_{ij}=0$ or $A_{ji}=0$. Since positive semidefinite matrices are Hermitian, $A$ is diagonal.
\end{proof}

By the construction of $G(A)$, if $A$ is a nonnegative matrix, we can see that $G(A)=G(\phi(A))$. Moreover, $\phi(A)-diag((\phi(A))_{11},\dots,(\phi(A))_{nn})$ is also the adjacency matrix of $G(A)$.
The following Theorem is a sufficient condition for $\operatorname{per}(A)=\lvert \det(A) \rvert$.
\begin{thm}\label{thm1} Let $A$ be a $\{0,1\}$-matrix. Suppose that the following statements hold:
	\begin{enumerate}
		\item There exists $\tau\in S_n$ such that $(AP_\tau)_{ii} = 1$ for all $i\in[n]$.
		\item For each $k = 1,\dots,\lfloor\frac{n}{2}\rfloor$ and each $i\in[n]$, we have $(B^{2k})_{ii}=((B^k)_{ii})^2$, where $B = AP_\tau-I$.
	\end{enumerate}
	Then $0\neq \operatorname{per}(A) = \operatorname{sgn}(\tau)\det(A)$.
\end{thm}
\begin{proof}
	Suppose that (1) and (2) occur. By (1), it is obvious that $\operatorname{per}(A)\neq 0$. Let $\theta\in S_n$ such that $L_\theta(A) = 1$. Since $A_{i \,\tau(i)}=(AP_\tau)_{i\,\tau^{-1}\theta(i)}$ for each $i\in[n]$, $$1=L_\tau(A) = L_{\tau^{-1}\theta}(AP_\tau).$$  Note that $B_{ii}=0$ for all $i\in[n]$. Let $\tau^{-1}\theta = c_1\cdots c_r$ be the disjoint cycles decomposition. If $c_1$ is a transposition, says $(jl)$, then $B_{jl}$ and $B_{lj}$ are positive. This means that $(B^2)_{jj} \geq B_{jl}B_{lj} > 0$. This cannot occur because of (2) with $k=1$. Thus $c_1$ is not a transposition; that is, $c_1$ is a cycle of length at least $3$. Then $G(AP_\tau)$ contains the dicycle that corresponding to $c_1$. We know that $(B^{2k})_{ii}$ is the number of walks of length $2k$ whose origin vertex and terminus vertex are $i$, while $((B^k)_{ii})^2$ is equal to the number of walks of length $2k$ whose origin vertex and terminus vertex are $i$ obtained by joining any two walks of length $k$ whose origin vertex and terminus vertex are $i$. By (2), $G(B)$ contains no cycle of even length. Thus, $c_1$ is a cycle of odd length. Since the composition of any two disjoint cycles is commutative, $c_1$ can be viewed arbitrarily; namely, $c_1,\dots,c_r$ are all cycles of odd length. This implies that $\tau^{-1}\theta$ is an even permutation. So, $\theta$ and $\tau$ have the same sign. Since $\theta$ is arbitrary for which $L_\theta(A)\neq 0$, $\operatorname{per}(A)=\operatorname{sgn}(\tau)\det(A)$.
\end{proof}

By Lemma \ref{cor22}, for a given complex matrix $A$, $\operatorname{per}(A) =  \operatorname{sgn}(\tau)\det(A)$ if $\phi(A)$ satisfies (1) and (2) in Theorem \ref{thm1}.

\begin{exaple}\label{ex1}
	Let
	$$\Gamma=\begin{pmatrix}
	0 & 0 & \ast & \ast & 0\\
	0 & 0 & 0 & \ast & \ast\\
	\ast & 0 & \ast & 0 & \ast\\
	\ast & \ast & 0 & 0 & 0\\
	0 & \ast & \ast & 0 & 0
	\end{pmatrix}.$$
	
	We have that
	$$\phi(\Gamma)=\begin{pmatrix}
	0 & 0 & 1 & 1 & 0 \\
	0 & 0 & 0 & 1 & 1  \\
	1 & 0 & 1 & 0 & 1 \\
	1 & 1 & 0 & 0 & 0 \\
	0 & 1 & 1 & 0 & 0 
	\end{pmatrix}.$$
	
	Choose $\tau = (13524)$, we can see that 
	$$B:=\phi(\Gamma)P_\tau-I = \begin{pmatrix}
	0 & 1 & 0 & 0 & 0 \\
	0 & 0 & 1 & 0 & 0  \\
	1 & 0 & 0 & 1 & 0 \\
	0 & 0 & 0 & 0 & 1 \\
	1 & 0 & 0 & 0 & 0 
	\end{pmatrix}.$$
	
	By calculation, we have
	$$B^2 = \begin{pmatrix}
	0 & 0 & 1 & 0 & 0 \\
	1 & 0 & 0 & 1 & 0  \\
	0 & 1 & 0 & 0 & 1 \\
	1 & 0 & 0 & 0 & 0 \\
	0 & 1 & 0 & 0 & 0 
	\end{pmatrix} \text{ and } B^4=\begin{pmatrix}
	0 & 1 & 0 & 0 & 1 \\
	1 & 0 & 1 & 0 & 0  \\
	1 & 1 & 0 & 1 & 0 \\
	0 & 0 & 1 & 0 & 0 \\
	1 & 0 & 0 & 1 & 0 
	\end{pmatrix}.$$
	Thus, $(B^{2k})_{ii}=((B^k)_{ii})^2$, where $k =1,2$. Because $\operatorname{sgn}(\tau)=1$, by Theorem \ref{thm1}, $\operatorname{per}(A)=\det(A)$. By Lemma \ref{cor22}, we also have that $\operatorname{per}(A)=\det(A)$ for all $A \in Q(\Gamma)$. 
\end{exaple}

Now, if $A$ is a nonnegative matrix such that $\operatorname{per}(A)=\det(A)$, we focus on the equality of $\operatorname{per}(A^k)$ and $\det(A^k)$. We know that $\det(A^k)=(\det(A))^k$. However, this property is not necessarily true for the permanent. So, $\operatorname{per}(A)=\det(A)$ does not mean $\operatorname{per}(A^k)=\det(A^k)$.

\begin{thm}\label{thm2}
	Let $A$ be an $n\times n$ nonnegative matrix such that $\operatorname{per}(A)=\det(A)$ and $A_{ii}\neq 0$ for all $i\in[n]$. Denote $e$ the identity element of $S_n$. Then the following statements are equivalent:
	\begin{enumerate}
		\item $\operatorname{per}(A^k)=\det(A^k)$ for all $k \in\mathbb{N}$.
		\item $\operatorname{per}(A^2)=\det(A^2)$.
		\item $\operatorname{per}(A)=L_e(A)$.
		\item $\operatorname{per}(A^k)=L_e(A^k)$ for all $k\in\mathbb{N}$.
	\end{enumerate}
\end{thm}
\begin{proof}
	It is obvious that (1) implies (2) and (4) implies (1).
	
	Suppose that (3) is not true. Then there exists a cycle $\tau\in S_n$ such that $L_\tau(A)\neq 0$. Since $\operatorname{per}(A)=\det(A)$, $\tau$ must be a cycle of odd length. So, there exists $i_0 \in [n]$ such that $\tau=(i_0 \, \,\tau(i_0)\, \tau^2(i_0)\, \dots\, \tau^{m-1}(i_0))$, where $m$ is the length of $\tau$. Let $D=diag(A_{11},\dots,A_{nn})$. Define an $n \times n$ matrix $B$ by
	\begin{equation*}
	B_{ij} = \left\{
	\begin{array}{ll}
	A_{ij}, & \hbox{if $i\neq j = \tau(i)$,} \\
	0, & \hbox{otherwise.}
	\end{array}
	\right.
	\end{equation*}
	Then $A$ can be considered as 
	$$A=D+B+X$$
	for some nonnegative matrix $X$. So, 
	$$A^2 = D^2 + DB + B^2 + Y$$ 
	for some nonnegative matrix $Y$. Note that $(B^2)_{i_0 \, \tau^2(i_0)} = B_{i_0 \, \tau(i_0)}B_{\tau(i_0)\,\tau^2(i_0)}\neq 0$. It follows that $(A^2)_{i_0 \, \tau^2(i_0)}\neq 0$. For each $i\in[n]$, since $(D^2)_{ii}\neq 0$, we have $(A^2)_{ii}\neq 0$. Similarly, because $B_{\tau^s(i_0) \,\tau^{s+1}(i_0)}\neq 0$ for all $s\in[m-2]$, $DB_{\tau^s(i_0) \,\tau^{s+1}(i_0)}$, and thus $A_{\tau^s(i_0) \,\tau^{s+1}(i_0)}$ are not zero. Let 
	$$\theta = (i_0 \tau^2(i_0) \tau^3(i_0) \dots \tau^{m-1}(i_0)).$$
	We deduce that $\theta$ is a cycle of even length such that $L_\theta(A^2)\neq 0$. So, $\operatorname{per}(A^2)\neq\det(A^2)$. Therefore, (2) implies (3).
	
	Now, suppose that (3) occurs. It is obvious that (4) is true if $A$ is a diagonal matrix. Suppose that $A$ is not diagonal; that is, $A_{ij}\neq 0$ for some distinct $i,j\in[n]$.  For each cycle $c$ in $S_n$, then $c\neq e$.  By the assumption (3) and the assumption that $A_{ii}>0$ for all $i\in [n]$, we conclude that $L_c(A)=0$; namely, there exists $i\in[n]$ such that $A_{i c(i)}=0$. Let $\tau=(i_0\, \tau(i_0)\, \dots\, \tau^{m-1}(i_0))$ be the cycle of length $m$ such that $A_{i_0\, \tau(i_0)}, A_{\tau(i_0)\,\tau^2(i_0)},\dots,A_{\tau^{m-2}(i_0)\,\tau^{m-1}(i_0)}$ are not zero. Now, we write $A = D + X$, where $D = diag(A_{11},\dots,A_{nn})$ and $X = A-D$. Thus, for each positive integer $k \geq 2$,
	$$A^k = (D+X)^k = D^k+D^{k-1}X + Y,$$
	for some nonnegative matrix $Y$. This implies that $(A^k)_{ii} \neq 0$ for all $i \in [n]$. Moreover,  $X_{i_0 \, \tau(i_0)}=A_{i_0\,\tau(i_0)}$ and $X_{\tau^s(i_0)\, \tau^{s-1}(i_0)} = A_{\tau^s(i_0)\, \tau^{s-1}(i_0)} \neq 0$ for all $s \in [m-2]$. These facts and the term $D^{k-1}X$ of $A^k$ implies that $$(A^k)_{i_0\, \tau(i_0)}, (A^k)_{\tau(i_0)\,\tau^2(i_0)},\dots,(A^k)_{\tau^{m-2}(i_0)\,\tau^{m-1}(i_0)}$$ are not zero. Suppose that there exists a positive integer $r \geq 2$ such that $(A^r)_{\tau^{m-1}(i_0), i_0}\neq 0$. There must exist $t_1\in[n]$ such that $(A^{r-1})_{\tau^{m-1}(i_0)\,t_1}$ and $A_{t_1 \,i_0}$ are not zero. By the similar reason, there exists $t_2\in[n]$ such that $(A^{r-2})_{\tau^{m-1}(i_0)\,t_2}$ and $A_{t_2\,t_1}$ are not zero. By repeating this process, there exists $t_1,t_2,\dots,t_{r-1}$ such that $A_{\tau^{m-1}(i_0)\,t_{r-1}}, A_{t_{r-1}\,t_{r-2}},\dots,A_{t_1\,i_0}$ are not zero. Thus, 
	$$i_0\,\tau(i_0)\,\dots\,\tau^{m-1}(i_0)\,t_{r-1}\,t_{r-2}\,\dots\,t_1\, i_0$$
	is a walk from $i_0$ to itself in $G(A)$. Since every walk from $u$ to $v$ contains a path from $u$ to $v$ and $A_{\tau(i_0)\,i_0}=0$ (by Theorem \ref{lem1}), there exists a cycle $\theta$ for which $\theta\neq e$ such that $L_\theta(A) \neq 0$. This contradicts to (3). Therefore $(A^k)_{\tau^{m-1}(i_0)\, i_0}= 0$ (i.e., $L_\tau(A^k)=0$) for all positive integer $k$.  Since $\tau$ is arbitrary cycle in $S_n$, $\sigma\neq e$ implies $L_\sigma(A^k)=0$ for all $k\in \mathbb{N}$.  Therefore (4) holds.
\end{proof}

Note that Theorem \ref{thm2} requires the assumption that $A_{ii}\neq 0$. The following theorem is the improvement for which a nonnegative matrix does not need that property.
\begin{thm}\label{thm3}
	Let $A$ be a nonnegative matrix such that $\operatorname{per}(A)\neq 0$. Then there exists $\tau \in S_n$ such that $L_e(AP_\tau)\neq 0$, where $e$ is the identity element of $S_n$. Then $\operatorname{per}(A^k)=\det(A^k)$ for all $k\in\mathbb{N}$ if and only if the following statements hold:
	\begin{enumerate}
		\item $\tau\in A_n$.
		\item $\operatorname{per}(A) = L_\tau(A)$.
		\item $\operatorname{per}(A^s)=L_e(A^s)$, where $s$ is the order of $\tau$.
	\end{enumerate}
\end{thm}
\begin{proof}
	Suppose that $\operatorname{per}(A^k)=\det(A^k)$ for all $k \in \mathbb{N}$. Because $\operatorname{per}(A)=\det(A)$, it is obvious that $\tau \in A_n$. If we consider $A$ as $DP_{\tau^{-1}} + X$ for some nonnegative matrix $X$ and invertible diagonal matrix $D$, we can see that
	$$A^s = (DP_{\tau^{-1}} + X)^s = (DP_{\tau^{-1}})^s+Y,$$
	for some nonnegative matrix $Y$. Note that $DP_{\tau^{-1}} = P_{\tau^{-1}}\acute{D}$ for some invertible diagonal matrix $\acute{D}$. Thus $(DP_{\tau^{-1}})^s = (P_{\tau^{-1}})^s\tilde{D} = \tilde{D}$, for some invertible diagonal matrix $\tilde{D}$. This implies that $(A^s)_{ii}\neq 0$ for each $i \in [n]$. Since $\operatorname{per}(A^k)=\det(A^k)$ for all $k\in\mathbb{N}$, we also have that $\operatorname{per}((A^s)^k)=\det((A^s)^k)$ for all $k \in \mathbb{N}$. By Theorem \ref{thm2}, (3) is true. To prove (2), suppose that there exists $\theta \in A_n$ such that $L_\theta(A)\neq 0$. Then there exists invertible diagonal matrices $D_1,D_2$ and a nonnegative matrix $X$ such that $$A=D_1P_{\tau^{-1}}+D_2P_{\theta^{-1}}+X.$$ Thus, because $$A^s=(D_1P_{\tau^{-1}}+D_2P_{\theta^{-1}}+X)^s = (D_1P_{\tau^{-1}})^{s-1}(D_2P_{\theta^{-1}})+Y= D_3P_{\tau\theta^{-1}}+Y$$ for some nonnegative matrix $Y$ and some invertible diagonal matrix $D_3$. Since (3) holds, $D_3P_{\tau\theta^{-1}}$ must be diagonal; that is, $P_{\tau\theta^{-1}}$ is the identity matrix. This happens if and only if $\theta=\tau$. Hence (2) holds.
	
	Conversely, suppose that (1)-(3) occur. By Theorem \ref{thm2} again, we have that $L_e((A^s)^m)=\operatorname{per}((A^s)^m)=\det((A^s)^m)$ for all $m\in\mathbb{N}$. This implies that
	\begin{equation}\label{eqpower}
	\hbox{ $\operatorname{per}(A^{ms})=\det(A^{ms})$ for all $m \in \mathbb{N}$.}	
	\end{equation}
	Suppose that $\operatorname{per}(A^k)\neq \det(A^k)$ for some $k \in \mathbb{N}$.  Since $\operatorname{per}(A^k)\neq \det(A^k)$, there is an odd permutation $\theta\in S_n$ such that $L_\theta(A^k)\neq 0$. We can consider $A^k$ as
	$$A^k=D_1P_{\theta^{-1}}+X_1$$
	for some invertible diagonal matrix $D_1$ and nonnegative matrix $X_1$. By (2), $A$ can be viewed as $A=D_2P_{\tau^{-1}}+X_2$ for some invertible diagonal matrix $D_2$ and nonnegative matrix $X_2$. By the division algorithm, $k=qs+r$, where $0\leq r <s$. We consider $$A^{s-r}=(D_2P_{\tau^{-1}}+X_2)^{s-r}=(D_2P_{\tau^{-1}})^{s-r}+X_3=D_3P_{(\tau^{-1})^{s-r}}+X_3$$ for some invertible diagonal matrix $D_3$ and nonnegative matrix $X_3$. Since $(q+1)s = k+s-r$, there also exists a invertible diagonal matrix $D_4$ and a nonnegative $X_4$ such that
	$$A^{(q+1)s}=A^kA^{s-r}=(D_1P_{\theta^{-1}}+X_1)(D_3P_{(\tau^{-1})^{s-r}}+X_3)=D_4P_{\theta^{-1}(\tau^{-1})^{s-r}}+X_4.$$
	This yields that $L_\sigma(A^{(q+1)s})\neq 0$, where $\sigma=\theta^{-1}(\tau^{-1})^{s-r}$.
	Since $\tau$ is even and $\theta$ is odd, we have $\sigma=\theta^{-1}(\tau^{-1})^{s-r}$ is an odd permutation.  Consequently, $\operatorname{per}(A^{(q+1)s})\neq \det(A^{(q+1)s})$, which is a contradiction to (\ref{eqpower}) . Hence (2) holds.
\end{proof}

\begin{exaple}\label{ex2}
	
	Let $$A=\begin{pmatrix}
	0 & 0 & 0 & 1 & 0 & 0 \\
	1 & 0 & 0 & 0 & 0 & 0 \\
	1 & 1 & 0 & 1 & 1 & 0 \\
	0 & 1 & 0 & 0 & 0 & 0 \\
	1 & 1 & 0 & 1 & 0 & 1 \\
	1 & 1 & 1 & 1 & 0 & 0
	\end{pmatrix}.$$ We can see that $\operatorname{per}(A)=\det(A)=L_\sigma(A)=1$, where $\sigma = (142)(356)$. By the calculation, we have
	$$A^{3}=\begin{pmatrix}
	1 & 0 & 0 & 0 & 0 & 0 \\
	0 & 1 & 0 & 0 & 0 & 0 \\
	3 & 3 & 1 & 3 & 0 & 0 \\
	0 & 0 & 0 & 1 & 0 & 0 \\
	3 & 3 & 0 & 3 & 1 & 0 \\
	3 & 3 & 0 & 3 & 0 & 1
	\end{pmatrix}.$$
	Since $\operatorname{per}(A^{3})=\det(A^{3})=L_e(A^3)=1$, by Theorem \ref{thm3}, $\operatorname{per}(A^k) = \det(A^k)$ for all $k \in \mathbb{N}$.
\end{exaple}

By Lemma \ref{cor22}, we can remark that, if $\operatorname{per}(A^k)=\det(A^k)$,  $\operatorname{per}(B^k)=\det(B^k)$ every matrix $B$ with the same support as $A$; that is, if $B$ is the matrix of the pattern
$$A=\begin{pmatrix}
0 & 0 & 0 & \ast & 0 & 0 \\
\ast & 0 & 0 & 0 & 0 & 0 \\
\ast & \ast & 0 & \ast & \ast & 0 \\
0 & \ast & 0 & 0 & 0 & 0 \\
\ast & \ast & 0 & \ast & 0 & \ast \\
\ast & \ast & \ast & \ast & 0 & 0
\end{pmatrix},$$
we can conclude that $\operatorname{per}(B^k)=\det(B^k)$ for every $k\in[n]$.

It is well-known that $\det(A^k)=(\det(A^k))$ for each square matrix $A$. So, if we assure that $A$ satisfies (1)-(3) in Theorem \ref{thm3}, the permanent of $A^k$ can be obtained from $\operatorname{per}(A)$ for each $k \in \mathbb{N}$. However, when the size of a matrix is larger, verifying those conditions also seem to be more difficult. The following theorem may help in this situation.

\begin{thm}
	Let $A$ be an $n\times n$ nonnegative matrix such that $\operatorname{per}(A)\neq 0$. Then $\operatorname{per}(A)=L_\tau(A)$ for some $\tau \in S_n$ if and only if there exists permutation matrices $P,Q$ such that $PAQ$ is upper triangular.
\end{thm}
\begin{proof}
	Suppose that $\operatorname{per}(A)=L_\tau(A)$ for some $\tau \in S_n$. If each row and each column of $A$ contain at least 2 nonzero entries, then there exists $\sigma \in S_n$ such that $\sigma \neq \tau$ and $L_\sigma(A)\neq 0$, a contradiction. Thus, $A$ must have a row or column that contains exactly one nonzero entry. If it is a row, choose permutation matrices $P_1,Q_1$ those permute that entry to $(P_1AQ_1)_{nn}$. In case it is a column, choose $P_1,Q_1$ so that the nonzero entry is permuted to $(P_1AQ_1)_{11}$. Now, consider the principal submatrix of $P_1AQ_1$ obtained by deleting that entry. Denote $B$ be this principal submatrix of $P_1AQ_1$. Since $\operatorname{per}(A)=L_\tau(A)$ for some $\tau \in S_n$, we have that $\operatorname{per}(B)=L_\theta(B)$ for some $\theta \in S_{n-1}$. So, $B$ also contains a row or column with exactly one nonzero entry. Choose $P_2,Q_2$ in the same way. By repeating this method, we derive permutation matrices $P,Q$ such that $PAQ$ is upper triangular. The converse statement of this theorem is obvious. Thus the proof is completed.
\end{proof}

By Example \ref{ex2}, we can see that, if $\tau=(16)(2453)$ and $\theta=(1364)(25)$, $P_\tau A P_\theta$ is upper triangular.

\section*{Acknowledgements}
The authors would like to thank anonymous referee(s) for reviewing this manuscript.

\end{document}